\documentclass[11pt,a4paper]{article}

\usepackage[T1]{fontenc}
\usepackage{amsmath,amsthm,mathtools}
\usepackage{newtxtext,newtxmath}
\usepackage{graphicx}
\usepackage{float}
\usepackage{geometry}
\usepackage{enumitem}
\usepackage{booktabs}
\usepackage{array}
\usepackage{microtype}
\usepackage{hyperref}
\hypersetup{colorlinks=true,linkcolor=black,citecolor=black,urlcolor=black,
  pdftitle={A novel Non-Homogeneous Local Representation of the Braid Group},
  pdfauthor={Arash Pourkia}}
\numberwithin{equation}{section}
\allowdisplaybreaks
\newtheorem{definition}{Definition}[section]
\newtheorem{remark}[definition]{Remark}
\newtheorem{proposition}[definition]{Proposition}

\newtheorem{theorem}[definition]{Theorem}

\title{A novel Non-Homogeneous Local Representation of the Braid Group}
\author{Arash Pourkia\\[3pt]
\normalsize College of Engineering and Technology, American University of the Middle East, Kuwait\\[2pt]
\normalsize \href{mailto:arash.pourkia@aum.edu.kw}{arash.pourkia@aum.edu.kw}}
\date{}

\begin{document}
\maketitle

\begin{abstract}
We construct a non-homogeneous local representation of the braid group $B_n$ by tracking the labels of individual strands through successive crossings. Each generator acts through a $2\times2$ block whose position is determined by the crossing and whose coefficients depend on the ordered labels of the participating strands. Starting with the ansatz
\[
M_{xy}=\begin{pmatrix}a_x&a_y\\ b_{xy}&c_{xy}\end{pmatrix},
\]
for $n\geq3$, the cubic braid relations for arbitrary ordered triples of strands determine $c_{xy}=0$ and $b_{xy}=(1-a_x)/a_x$ in the nondegenerate branch. The moving-label rule leads to an explicit semidirect-product representation, for which we prove the Artin relations and use exact evaluations of the resulting symbolic matrices. In the special homogeneous specialization, where $a_x=a$ for every strand label $x$, the resulting matrix representation is diagonally equivalent to unreduced Burau. We prove that a non-homogeneous specialization of our representation detects Bigelow's explicit nontrivial $B_5$ Burau-kernel commutator. Thus the strand-dependent coefficients distinguish a nontrivial braid that is invisible to the classical Burau representation. This separation provides evidence that the non-homogeneous strand-tracking construction is a promising candidate for a faithful Burau-type representation.
\end{abstract}


\section{Introduction}\label{sec:introduction}

The braid group $B_n$ admits several familiar descriptions: it is the group of isotopy classes of braids on $n$ strands, the fundamental group of an unordered configuration space, and the mapping class group of an $n$-punctured disk. Artin's presentation relates these pictures to generators and relations \cite{Artin1925,Artin1947}; standard references include \cite{Birman1974,KasselTuraev2008}. The search for matrix representations that retain useful information about braids has been a recurring problem in this subject.

The Burau representation \cite{Burau1935} is a natural point of comparison. Long and Paton established its nonfaithfulness for $n\geq6$ \cite{LongPaton1993}, and Bigelow exhibited an explicit nontrivial element of its kernel for $n=5$ \cite{Bigelow1999}. In a 2026 preprint, Bharathram, Birman, and Brendle announced a proof of faithfulness in the formerly unresolved case $n=4$ \cite{BharathramBirmanBrendle2026}. Independently, the Lawrence--Krammer representation is faithful for every braid group \cite{Lawrence1990,Bigelow2001,Krammer2002}. Recent investigations of Burau and related braid-group representations include complex specializations and $q$-rationals for $B_3$ \cite{MorierGenoudOvsienkoVeselov2024}, kernel elements for evaluated Gassner representations \cite{Nasser2025Gassner}, irreducibility and equivalence of local matrix families \cite{ChreifDally2024,Nasser2024Irreducibility}, and extensions to singular and welded braids \cite{BardakovChbiliKozlovskaya2024,Sysoeva2024Welded,Nasser2026LocalExtensions}. Together with recent work on charge-conserving braid representations \cite{MartinRowell2026}, these results illustrate the continuing interest in braid-group representations and their connections to other areas of mathematics.

Braid representations also enter the study of knot and link invariants, with the Markov moves providing the connection to braid closures. Kauffman's \emph{Knots and Physics} \cite{Kauffman2012KnotsPhysics} and Adams's \emph{The Knot Book} \cite{Adams2004KnotBook} provide classical accounts of this relationship. The present work develops a strand-sensitive representation of the braid group and studies its power to distinguish braid elements.

The starting point of this paper is the path-based construction in \cite{Pourkia2022}; related approaches to local braid representations, from different points of view, appear in \cite{Mikhalchishina2013,TongYangMa1996}. In \cite{Pourkia2022}, the matrix assigned to a braid records information obtained by following its strands. In particular, the homogeneous $2\times2$ blocks satisfying the ordinary local braid relations are classified within the setting of that work, and a specialization yields unreduced Burau. In the present work, the local block continues to act on the two positions involved in a crossing, but its entries are determined by the labels of the strands currently occupying those positions.

Suppose that the crossing at Artin position $i$ involves strands labelled $x$ and $y$, in that order. We write $\widehat{\sigma}_{i;xy}$ for the associated local matrix. The index $i$ determines the position of its $2\times2$ block in an $n\times n$ matrix; the indices $x,y$ determine the coefficients. As the word is read, each crossing changes the ordering of the strands. The corresponding matrix factors must therefore be multiplied with the accompanying relabelling. We express this rule as a semidirect product of a matrix group by the symmetric group $S_n$.

For the local block we prescribe the top row to be $(a_x,a_y)$ and leave the lower entries undetermined:
\[
M_{xy}=\begin{pmatrix}a_x&a_y\\b_{xy}&c_{xy}\end{pmatrix}.
\]
This is our starting ansatz: the top row of the $2\times2$ local block records the two strands participating in the crossing. For $n\geq3$, invertibility and the labelled braid relations force $c_{xy}=0$ and $b_{xy}=(1-a_x)/a_x$ in the nondegenerate branch of this ansatz. We derive these relations explicitly, construct the corresponding semidirect-product representation, and verify the Artin relations.

Our principal test is Bigelow's explicit commutator in $B_5$, a nontrivial braid in the kernel of Burau. At $a_1=a_2=a_3=a_4=2$ and $a_5=-1$, the strand-tracking representation assigns this braid a nonidentity matrix, in contrast with the special homogeneous specialization $a_x=a$ for every strand label $x$, which is Burau-equivalent. This exact separation result demonstrates that allowing the coefficients to depend on the moving strands can detect a braid in the classical Burau kernel.

We establish the notation in Section~\ref{sec:preliminaries} and discuss the relevant literature in Section~\ref{sec:literature}. In Section~\ref{sec:main-results} we derive the labelled local block directly from the cubic relation for an arbitrary triple and general $n$. Section~\ref{sec:formal} defines the semidirect-product representation and states the resulting braid-group homomorphism. Section~\ref{sec:examples} specializes the local matrices to $B_5$ and gives the exact proof of Bigelow separation; Section~\ref{sec:conclusion} concludes.

\section{Preliminaries and notation}\label{sec:preliminaries}

For $n\geq2$, the braid group $B_n$ is generated by
\[
\sigma_1,\ldots,\sigma_{n-1}
\]
subject to the relations
\begin{equation}\label{eq:far}
\sigma_i\sigma_j=\sigma_j\sigma_i\qquad (|i-j|\geq2)
\end{equation}
and
\begin{equation}\label{eq:cubic}
\sigma_i\sigma_{i+1}\sigma_i
=
\sigma_{i+1}\sigma_i\sigma_{i+1}
\qquad(1\leq i\leq n-2).
\end{equation}
Throughout, braid diagrams are oriented downward and are read from top to bottom.  If
\[
\beta=\sigma_{i_1}^{\varepsilon_1}\sigma_{i_2}^{\varepsilon_2}\cdots\sigma_{i_m}^{\varepsilon_m},
\qquad \varepsilon_r\in\{+1,-1\},
\]
then the first displayed letter is the topmost crossing, the second displayed letter is the next crossing below it, and so on.  Matrix products associated with a braid word are written in this same left-to-right order.

We write
\[
[n]=\{1,2,\ldots,n\}
\]
for the strand-label set.  The strands beginning at the top endpoints are denoted
\[
s_1,s_2,\ldots,s_n,
\]
from left to right, and each label travels with its strand throughout the braid.  If, at some horizontal level, the strands occupying the positions from left to right are
\[
s_{x_1},s_{x_2},\ldots,s_{x_n},
\]
we call
\[
(x_1,x_2,\ldots,x_n)
\]
the \emph{current strand order}.  In particular, a local order $(x,y,z)$ means that the three adjacent positions under discussion are occupied, from left to right, by $s_x,s_y,s_z$.  The initial top order is $(1,2,\ldots,n)$.

Let $S_n$ be the symmetric group.  For $1\leq i\leq n-1$, let
\[
P_i=(i\;i+1)\in S_n.
\]
The identity permutation is denoted by $\mathrm{id}$.  Permutation multiplication is composed according to
\[
(PQ)(r)=P(Q(r)).
\]
Thus, when a product $PQ$ acts on a label, $Q$ acts first and then $P$.  The elementary transpositions satisfy the Coxeter relations
\begin{equation}\label{eq:coxeter-Sn}
P_i^2=\mathrm{id},\qquad
P_iP_j=P_jP_i\quad(|i-j|\geq2),\qquad
P_iP_{i+1}P_i=P_{i+1}P_iP_{i+1}.
\end{equation}
The last identity is the same three-term braid relation as in $B_n$, now inside the symmetric group.

For a braid or partial braid, its current strand order can equivalently be encoded by the permutation $P\in S_n$ whose value $P(r)$ is the label of the strand occupying position $r$.  Thus the top order corresponds to $\mathrm{id}$, and a single crossing $\sigma_i$ changes the order by the elementary transposition $P_i$.

The distinction between \emph{positions} and \emph{strand labels} is fundamental.  The Artin index $i$ records the adjacent positions in which a crossing occurs.  The labels $x,y$ record the actual strands occupying those positions at that stage of the braid.  We write
\[
\widehat{\sigma}_{i;xy}
\]
for the matrix assigned to a positive crossing at position $i$ at which strand $s_x$ passes over strand $s_y$.  The hat distinguishes these matrices from the abstract Artin generators $\sigma_i$.

For a field $F$, $M_n(F)$ denotes the algebra of $n\times n$ matrices over $F$, $GL_n(F)$ its group of invertible matrices, and $I_n$ the $n\times n$ identity matrix.

\section{Related literature and context}\label{sec:literature}

Recent work on braid-group representations, Burau and Gassner specializations, questions of faithfulness, and connections with knot theory shows that these subjects remain active areas of research. The developments range from the classification of particular families of matrix representations to braid actions on free groups, generalized braid groups, and Alexander-type knot invariants. We first review recent results and then identify the classical constructions most closely related to the strand-labelled matrices studied here.

Recent classification results include the work of Martin and Rowell \cite{MartinRowell2026}, who classified charge-conserving spin-chain braid representations within a specified monoidal framework. Martin, Rowell, and Torzewska \cite{MartinRowellTorzewska2025} established a related classification for charge-conserving loop braid representations. Hietarinta, Martin, and Rowell \cite{HietarintaMartinRowell2024} studied solutions of the constant Yang--Baxter equation under an additive charge-conservation ansatz in dimension three. These are systematic classifications within tensor or monoidal settings, distinct from the present $n\times n$ strand-labelled local construction. Another construction method, based on iterated twisted tensor products of algebras, was developed by Gustafson, Kimball, Rowell, and Zhang \cite{GustafsonKimballRowellZhang2020}.

For Burau theory itself, Morier-Genoud, Ovsienko, and Veselov \cite{MorierGenoudOvsienkoVeselov2024} related the complex specializations of the $B_3$ Burau representation to $q$-deformed rational numbers and obtained faithfulness results for specified regions of the parameter plane. Bharathram, Birman, and Brendle \cite{BharathramBirmanBrendle2026} announced in a 2026 preprint that the Burau representation of $B_4$ is faithful. In the multivariable setting, Nasser \cite{Nasser2025Gassner} produced nontrivial kernel elements for evaluated Gassner representations at roots of unity, distinguishing these specializations from the generic representation. These results illustrate the importance of both the choice of coefficients and the group action on them.

Recent work on local and low-dimensional matrix families provides further comparisons. Chreif and Dally \cite{ChreifDally2024} investigated the irreducibility and equivalence of homogeneous local representations, including relations with Burau and the standard representation. Nasser \cite{Nasser2024Irreducibility} determined irreducibility conditions for a composition factor of a braid-group representation introduced by Bardakov and Bellingeri. Kahil and Abdulrahim \cite{KahilAbdulrahim2022} studied the irreducibility of representations originating with Egea and Galina and constructed further multiparameter families. Sysoeva \cite{Sysoeva2021} proved that $B_n$ has no irreducible complex representations of dimension $n+1$ for $n\geq10$; in a later preprint \cite{Sysoeva2024Welded}, she studied extensions of certain irreducible braid representations, including the reduced Burau and standard types, to welded braid groups. These results address representation-theoretic questions in settings different from the relabelling mechanism of our construction.

Related investigations enlarge the class of allowed crossings or braids. Bardakov, Chbili, and Kozlovskaya \cite{BardakovChbiliKozlovskaya2024} studied extensions of braid representations to singular braid monoids and groups, including an extension of the Lawrence--Krammer--Bigelow representation. Nasser \cite{Nasser2026LocalExtensions} compared local and $\Phi$-type extensions of specified homogeneous local representations to the singular braid monoid. Bardakov and Nasybullov \cite{BardakovNasybullov2020} developed multi-switch methods for virtual braid representations and virtual-link invariants. Bardakov and collaborators \cite{BardakovEtAl2022} constructed representations of braid and virtual braid groups related to the Lawrence--Krammer--Bigelow representation and obtained faithful representations of associated crystallographic quotients; these quotient statements are distinct from faithfulness of a representation of the full braid group.

The interplay between braid actions and knot invariants is equally active. G\"ug\"umc\"u and Kauffman \cite{GugumcuKauffman2025} constructed mock Alexander polynomials for starred links and linkoids in surfaces through diagrammatic state sums. Juh\'asz, Kauffman, and Ogasa \cite{JuhaszKauffmanOgasa2024} used spanning surfaces in thickened surfaces to define Seifert-matrix and Alexander-type invariants, including a polynomial relevant to virtual knots. Kauffman's \emph{Knots and Physics} \cite{Kauffman2012KnotsPhysics} provides the classical background in skein theory, state sums, and polynomial knot invariants. These approaches construct invariants of knots and diagrams; the present work instead uses matrix coordinates to distinguish braid elements.

Adams's \emph{The Knot Book} \cite{Adams2004KnotBook} presents the basic relationship among braids, closures, and knot invariants. More recently, Adams and collaborators \cite{AdamsEtAl2024GeneralisedKnotoids} developed generalised knotoids involving several poles, intervals, and circles, and considered analogues of index and bracket polynomials. Their constructions enlarge the types of knot-like diagrams under study, whereas our construction records the changing labels of strands within ordinary braids.

Free-group actions provide another connection between braid groups and knot theory. Shpilrain \cite{Shpilrain2024KnotInvariants} studied knot and link invariants obtained from the Artin action on free groups and from Wada's braid actions, comparing the resulting invariants with Alexander-polynomial constructions. The use of free-group automorphisms and Fox derivatives there contrasts with the direct strand-labelled matrices considered here, while providing a useful point of comparison for the information encoded by a braid action.

The foundational matrix models remain important. Burau's representation \cite{Burau1935} is not faithful for $n\geq5$, as established by the explicit examples of Long and Paton \cite{LongPaton1993} and Bigelow \cite{Bigelow1999}. The Lawrence--Krammer representation, by contrast, provides faithful linear representations of braid groups \cite{Lawrence1990,Bigelow2001,Krammer2002}. Formanek's low-degree analysis \cite{Formanek1996} and the Tong--Yang--Ma family \cite{TongYangMa1996} supply further classical comparisons.

Mikhalchishina \cite{Mikhalchishina2013} investigated local linear representations of $B_3$ and homogeneous local representations of $B_n$, including links with the Burau representation and constructions arising from Wada's free-group actions. This work provides a systematic benchmark when the same local block is used for every crossing position. The coefficients in our construction, by contrast, are changed according to the current strand order; the cited homogeneous-local analysis does not itself amount to a classification of such strand-dependent semidirect-product representations.

The closest historical precedents for strand-dependent parameters are the multivariable Gassner representation on pure braids \cite{Gassner1961} and the labelled multivariable Burau matrices studied by Penne \cite{Penne1995,Penne1998}. In these constructions, the coefficients associated with individual strands are transported through crossings. Our construction begins with a specified labelled $2\times2$ ansatz, derives its coefficients from the labelled braid relations, and formulates the strand permutation as part of an explicit semidirect-product representation. Nasser and Abdulrahim's work \cite{NasserAbdulrahim2021} on irreducibility of representations associated with extensions of Burau- and Gassner-type constructions offers a related, though distinct, representation-theoretic comparison.

Finally, the path-based homogeneous family of \cite{Pourkia2022} motivates the local-block presentation used here. Related homogeneous-local classifications and braid matrices are discussed in \cite{Mikhalchishina2013,TongYangMa1996}. In the present paper the $2\times2$ block remains attached to the crossing position, while its coefficients depend on the moving strand labels. In the special homogeneous specialization $a_x=a$ for every strand label $x$, the matrix representation is diagonally equivalent to unreduced Burau; the non-homogeneous Bigelow calculation gives a second, direct comparison.

\section{Main results: the strand-labelled local matrices}\label{sec:main-results}

\subsection{The local-block ansatz}\label{subsec:informal-motivation}

We begin with the usual block placement for a homogeneous $2$-local representation.  In the setting of Theorem~3.1 of \cite{Pourkia2022}, one assigns to every Artin generator the same $2\times2$ core, placed in the two rows and columns corresponding to the crossing position.  Here we retain the moving block but allow its entries to record which strands actually meet.

At the outset, a completely unrestricted $2\times2$ matrix would have four unrelated entries, for example
\[
\begin{pmatrix}
u_{xy}&v_{xy}\\ w_{xy}&z_{xy}\end{pmatrix}.
\]
In the strand-labelled ansatz, we set the top row equal to the two strand parameters $a_x,a_y$ and initially leave the lower row unrestricted.  Thus, when the current local order at positions $i,i+1$ is $(s_x,s_y)$ and $\sigma_i$ is positive, our initial core is
\begin{equation}\label{eq:initial-general-core}
M_{xy}
=
\begin{pmatrix}
a_x&a_y\\
b_{xy}&c_{xy}
\end{pmatrix},
\qquad x\neq y,
\end{equation}
where $b_{xy}$ and $c_{xy}$ are, at this stage, independent ordered-label parameters.  Its determinant is
\[
\det M_{xy}=a_xc_{xy}-a_yb_{xy}.
\]
We work in the nondegenerate branch, where all $a_r$ are nonzero and every local core is invertible.  No nonzero assumption is imposed on $c_{xy}$; the cubic relation will determine it.

We first show that the lower-right entry in \eqref{eq:initial-general-core} cannot survive the dynamically labelled cubic relation.  Let $x,y,z$ be pairwise distinct.  On the active three positions, the left side of the cubic relation is
\begin{align}
&\begin{pmatrix}
a_x&a_y&0\\ b_{xy}&c_{xy}&0\\0&0&1
\end{pmatrix}
\begin{pmatrix}
1&0&0\\0&a_x&a_z\\0&b_{xz}&c_{xz}
\end{pmatrix}
\begin{pmatrix}
a_y&a_z&0\\b_{yz}&c_{yz}&0\\0&0&1
\end{pmatrix}
\notag\\[1mm]
&\qquad=
\begin{pmatrix}
a_xa_y(1+b_{yz})&a_xa_yc_{yz}+a_xa_z&a_ya_z\\
a_xb_{yz}c_{xy}+a_yb_{xy}&a_xc_{xy}c_{yz}+a_zb_{xy}&a_zc_{xy}\\
b_{xz}b_{yz}&b_{xz}c_{yz}&c_{xz}
\end{pmatrix},
\label{eq:general-c-left}
\end{align}
whereas the right side is
\begin{align}
&\begin{pmatrix}
1&0&0\\0&a_y&a_z\\0&b_{yz}&c_{yz}
\end{pmatrix}
\begin{pmatrix}
a_x&a_z&0\\b_{xz}&c_{xz}&0\\0&0&1
\end{pmatrix}
\begin{pmatrix}
1&0&0\\0&a_x&a_y\\0&b_{xy}&c_{xy}
\end{pmatrix}
\notag\\[1mm]
&\qquad=
\begin{pmatrix}
a_x&a_xa_z&a_ya_z\\
a_yb_{xz}&a_xa_yc_{xz}+a_zb_{xy}&a_y^2c_{xz}+a_zc_{xy}\\
b_{xz}b_{yz}&a_xb_{yz}c_{xz}+b_{xy}c_{yz}&a_yb_{yz}c_{xz}+c_{xy}c_{yz}
\end{pmatrix}.
\label{eq:general-c-right}
\end{align}
Subtracting the right side from the left side gives
\begin{equation}\label{eq:general-c-difference}
\begin{pmatrix}
a_x\bigl(a_y(1+b_{yz})-1\bigr)&a_xa_yc_{yz}&0\\[1mm]
a_xb_{yz}c_{xy}+a_y(b_{xy}-b_{xz})
&a_x(c_{xy}c_{yz}-a_yc_{xz})&-a_y^2c_{xz}\\[1mm]
0&(b_{xz}-b_{xy})c_{yz}-a_xb_{yz}c_{xz}
&c_{xz}-a_yb_{yz}c_{xz}-c_{xy}c_{yz}
\end{pmatrix}.
\end{equation}
Two entries already determine the lower-right coefficients.  The $(1,2)$ and $(2,3)$ entries of \eqref{eq:general-c-difference} are
\[
a_xa_yc_{yz},
\qquad
-a_y^2c_{xz}.
\]
Since $a_x$ and $a_y$ are nonzero, the cubic relation forces
\[
c_{yz}=0,
\qquad
c_{xz}=0.
\]
Because the ordered triple $(x,y,z)$ may be chosen arbitrarily among pairwise distinct labels, this gives
\begin{equation}\label{eq:c-forced-zero}
\boxed{
c_{xy}=0
\qquad(x\neq y)
}
\end{equation}
for every ordered pair when $n\geq3$.  Indeed, for any $u\neq v$, choose a third label $w$ and apply the $(1,2)$ entry above to the ordered triple $(w,u,v)$.  For $B_2$ there is no cubic relation; as with the later restriction on $b_{xy}$, we adopt \eqref{eq:c-forced-zero} there as the uniform member of the same family.

Thus the zero in the lower-right corner is not assumed.  It is forced by the cubic braid relation in the nondegenerate strand-labelled branch.  From this point onward the preliminary core is therefore
\begin{equation}\label{eq:precore}
\boxed{
M_{xy}
=
\begin{pmatrix}
a_x&a_y\\
b_{xy}&0
\end{pmatrix},
\qquad x\neq y.
}
\end{equation}
The symbols $b_{xy}$ remain independent ordered-label parameters at this stage.  Proposition~\ref{prop:ordered-triple} below will show that requiring the cubic relation for all ordered triples determines these coefficients completely.  The corresponding $n\times n$ matrix is
\begin{equation}\label{eq:prehat}
\boxed{
\widehat{\sigma}_{i;xy}
=
I_{i-1}\oplus
\begin{pmatrix}
a_x&a_y\\
b_{xy}&0
\end{pmatrix}
\oplus I_{n-i-1}.
}
\end{equation}
Thus the block position records $i$, whereas its entries record strand labels $x$ and $y$.  Its determinant is now
\begin{equation}\label{eq:precore-det}
\det M_{xy}=-a_yb_{xy}.
\end{equation}
Consequently, invertibility in this branch requires $a_y\neq0$ and $b_{xy}\neq0$.

\begin{figure}[H]
\centering
\begin{minipage}[c]{0.28\textwidth}
\centering
\includegraphics{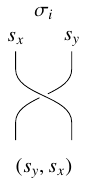}
\end{minipage}%
\hfill
\begin{minipage}[c]{0.08\textwidth}
\centering
$\longleftrightarrow$
\end{minipage}%
\hfill
\begin{minipage}[c]{0.58\textwidth}
\centering
$\displaystyle
\widehat{\sigma}_{i;xy}
=
I_{i-1}\oplus
\begin{pmatrix}
a_x&a_y\\
b_{xy}&0
\end{pmatrix}
\oplus I_{n-i-1}$
\end{minipage}
\caption{A positive crossing at Artin position $i$ with strand $s_x$ passing over strand $s_y$.  The block location records the crossing position, while the entries record the ordered strand labels.}
\label{fig:local-block}
\end{figure}

Figure~\ref{fig:local-block} summarizes the two pieces of information carried by a positive local matrix: the block position records the Artin index, while the entries record the ordered labels of the strands that actually interact.

\begin{figure}[H]
\centering
\begin{minipage}[c]{0.28\textwidth}
\centering
\includegraphics{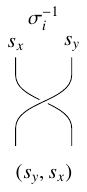}
\end{minipage}%
\hfill
\begin{minipage}[c]{0.08\textwidth}
\centering
$\longleftrightarrow$
\end{minipage}%
\hfill
\begin{minipage}[c]{0.60\textwidth}
\centering
$\displaystyle
\widehat{\sigma}_{i;yx}^{-1}
=
I_{i-1}\oplus
\begin{pmatrix}
0&b_{yx}^{-1}\\[1mm]
a_x^{-1}&-\dfrac{a_y}{a_xb_{yx}}
\end{pmatrix}
\oplus I_{n-i-1}$
\end{minipage}
\caption{A negative crossing $\sigma_i^{-1}$ when the current local order is $(s_x,s_y)$.  The inverse uses the swapped ordered labels $y,x$.  This is previewed here at the matrix level; after the semidirect-product bullet multiplication is defined, the inverse formula \eqref{eq:H-inverse} makes this relabelling explicit.}
\label{fig:negative-local-block}
\end{figure}

The matrix data of successive crossings will eventually be multiplied by ordinary matrix multiplication.  Because the strand labels change after each crossing, however, the multiplication law must also retain the natural strand permutation.  Figure~\ref{fig:negative-local-block} anticipates why this relabelling is essential for inverse letters.  The semidirect-product structure is formalized in Section~\ref{sec:formal}.  We next determine the remaining coefficient restriction required by the cubic relation.

\subsection{The general ordered triple \texorpdfstring{$(x,y,z)$}{(x,y,z)}}\label{subsec:general-triple}

Let $n\geq3$, let $1\leq i\leq n-2$, and suppose that the three adjacent positions $i,i+1,i+2$ are occupied, from left to right, by strands with pairwise distinct labels $(x,y,z)$. Along the two sides of the cubic relation the labels follow the successive orders
\[
\begin{aligned}
(x,y,z)&\xrightarrow{\sigma_i}(y,x,z)
\xrightarrow{\sigma_{i+1}}(y,z,x)
\xrightarrow{\sigma_i}(z,y,x),\\[2pt]
(x,y,z)&\xrightarrow{\sigma_{i+1}}(x,z,y)
\xrightarrow{\sigma_i}(z,x,y)
\xrightarrow{\sigma_{i+1}}(z,y,x).
\end{aligned}
\]
Thus the dynamically labelled matrix products are
\begin{align}
L_{xyz}^{(i)}
&=\widehat{\sigma}_{i;xy}\widehat{\sigma}_{i+1;xz}\widehat{\sigma}_{i;yz},
\label{eq:general-cubic-left-word}\\
R_{xyz}^{(i)}
&=\widehat{\sigma}_{i+1;yz}\widehat{\sigma}_{i;xz}\widehat{\sigma}_{i+1;xy}.
\label{eq:general-cubic-right-word}
\end{align}

\begin{figure}[H]
\centering
\includegraphics[width=0.91\linewidth]{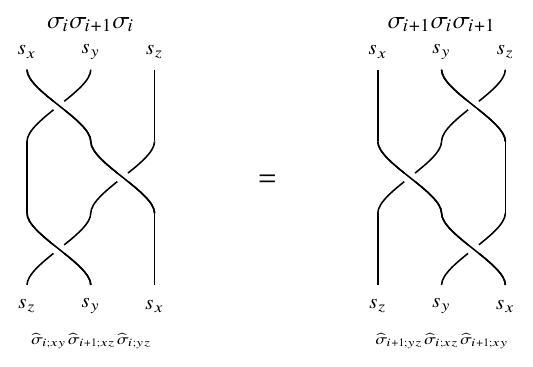}
\caption{The cubic relation at positions $i,i+1,i+2$ for the general initial strand order $(x,y,z)$. The two diagrams end in the same order $(z,y,x)$, but the local matrices must follow the successive ordered crossing pairs $(x,y),(x,z),(y,z)$ and $(y,z),(x,z),(x,y)$, respectively.}
\label{fig:cubic-general}
\end{figure}

By \eqref{eq:precore}, only the three active positions contribute nonidentity blocks. Ordinary multiplication gives
\begin{align}
L_{xyz}^{(i)}&=I_{i-1}\oplus
\begin{pmatrix}
a_xa_y(1+b_{yz})&a_xa_z&a_ya_z\\
a_yb_{xy}&a_zb_{xy}&0\\
b_{xz}b_{yz}&0&0
\end{pmatrix}
\oplus I_{n-i-2},\label{eq:general-left-matrix}\\[2mm]
R_{xyz}^{(i)}&=I_{i-1}\oplus
\begin{pmatrix}
a_x&a_xa_z&a_ya_z\\
a_yb_{xz}&a_zb_{xy}&0\\
b_{xz}b_{yz}&0&0
\end{pmatrix}
\oplus I_{n-i-2}.\label{eq:general-right-matrix}
\end{align}
Consequently,
\begin{equation}\label{eq:general-difference}
L_{xyz}^{(i)}-R_{xyz}^{(i)}
=0_{i-1}\oplus
\begin{pmatrix}
a_x\bigl(a_y(1+b_{yz})-1\bigr)&0&0\\
a_y(b_{xy}-b_{xz})&0&0\\
0&0&0
\end{pmatrix}
\oplus 0_{n-i-2},
\end{equation}
where a zero block of size $0$ is omitted at the boundary. The computation is therefore valid for arbitrary $n\geq3$, not only for a particular three- or five-strand example.

\begin{proposition}[Cubic criterion for an ordered triple]\label{prop:ordered-triple}
Let $n\geq3$, $1\leq i\leq n-2$, and let $x,y,z\in[n]$ be pairwise distinct. Under the nonzero-parameter convention above, the dynamically labelled cubic identity
\begin{equation}\label{eq:dynamic-cubic}
\widehat{\sigma}_{i;xy}
\widehat{\sigma}_{i+1;xz}
\widehat{\sigma}_{i;yz}
=
\widehat{\sigma}_{i+1;yz}
\widehat{\sigma}_{i;xz}
\widehat{\sigma}_{i+1;xy}
\end{equation}
holds if and only if
\begin{equation}\label{eq:ordered-triple-conditions}
\boxed{b_{xy}=b_{xz},\qquad b_{yz}=\frac{1-a_y}{a_y}.}
\end{equation}
\end{proposition}

\begin{proof}
The difference between the two sides is given by \eqref{eq:general-difference}. Since $a_x$ and $a_y$ are nonzero, the difference vanishes exactly when
\[
a_y(1+b_{yz})=1,\qquad b_{xy}=b_{xz},
\]
which is equivalent to \eqref{eq:ordered-triple-conditions}.
\end{proof}

The same criterion applies to every ordered triple of pairwise distinct labels. For $n\geq3$, requiring the cubic relation for every such triple forces
\begin{equation}\label{eq:b-restriction}
\boxed{b_{xy}=\frac{1-a_x}{a_x}\qquad(x,y\in[n],\ x\neq y).}
\end{equation}
Indeed, for any distinct labels $x,y$, choose $z\notin\{x,y\}$ and apply Proposition~\ref{prop:ordered-triple} to the ordered triple $(z,x,y)$. Its second condition is precisely \eqref{eq:b-restriction}. When $n=2$, there is no cubic relation, and \eqref{eq:b-restriction} is adopted as the uniform specialization defining the same family.

\subsection{The final two-parameter local block}\label{subsec:final-block}

After imposing the cubic restriction \eqref{eq:b-restriction}, the independent symbols $b_{xy}$ disappear from the coefficient system.  The local matrix attached to a crossing of strand $x$ over strand $y$ becomes
\begin{equation}\label{eq:final-core}
\boxed{
M_{xy}
=
\begin{pmatrix}
a_x&a_y\\[1mm]
\dfrac{1-a_x}{a_x}&0
\end{pmatrix}.
}
\end{equation}
Thus each local core depends on the two strand parameters $a_x$ and $a_y$.

Accordingly,
\begin{equation}\label{eq:final-hat}
\boxed{
\widehat{\sigma}_{i;xy}
=
I_{i-1}\oplus
\begin{pmatrix}
a_x&a_y\\[1mm]
\dfrac{1-a_x}{a_x}&0
\end{pmatrix}
\oplus I_{n-i-1}.
}
\end{equation}
Its determinant is
\begin{equation}\label{eq:final-det}
\det\widehat{\sigma}_{i;xy}
=
-a_y\frac{1-a_x}{a_x}.
\end{equation}
Over a numerical field specialization, invertibility therefore requires $a_x\neq0,1$ and $a_y\neq0$.  Since every label can occur in either role, it is natural to require
\[
a_r\neq0,1\qquad(r=1,\ldots,n).
\]
In the formal construction below we work symbolically over the rational function field $k(a_1,\ldots,a_n)$, where the determinant in \eqref{eq:final-det} is automatically nonzero.

\begin{remark}[Relation with the 2022 local family and unreduced Burau]\label{rem:burau-connection}
The homogeneous Case~1 family in Theorem~3.1 of \cite{Pourkia2022} has local core
\begin{equation}\label{eq:2022-case1}
C(a,b)=
\begin{pmatrix}
a&b\\[1mm]
\dfrac{1-a}{b}&0
\end{pmatrix}.
\end{equation}
The present strand-labelled construction does not recover the whole two-parameter family \eqref{eq:2022-case1} when the strand parameters are made homogeneous.  Indeed, if
\[
a_1=\cdots=a_n=a,
\]
then every local core \eqref{eq:final-core} becomes
\begin{equation}\label{eq:new-hom-core}
C_{\mathrm{hom}}(a)=
\begin{pmatrix}
a&a\\[1mm]
\dfrac{1-a}{a}&0
\end{pmatrix},
\end{equation}
which is precisely the specialization $b=a$ of \eqref{eq:2022-case1}.  The distinction is important: in the present construction the loss of the independent second parameter is a consequence of the dynamically labelled cubic relations, rather than an assumption made at the outset.

Nevertheless, the homogeneous family \eqref{eq:new-hom-core} recovers the full generic one-parameter unreduced Burau representation up to a diagonal change of basis.  Write
\[
t=1-a,
\qquad
B(t)=
\begin{pmatrix}
1-t&t\\
1&0
\end{pmatrix},
\]
so that $B(t)$ is the standard local block of the unreduced Burau representation.  Set
\[
r=\frac{a}{1-a}=\frac{1-t}{t},
\qquad
D_n=\operatorname{diag}(r^{n-1},r^{n-2},\ldots,r,1).
\]
For the adjacent entries of $D_n$ one has $d_i/d_{i+1}=r$.  Hence, for every $i$,
\begin{align}
D_n\bigl(I_{i-1}\oplus B(1-a)\oplus I_{n-i-1}\bigr)D_n^{-1}
&=I_{i-1}\oplus
\begin{pmatrix}
a&r(1-a)\\[1mm]
r^{-1}&0
\end{pmatrix}
\oplus I_{n-i-1}\notag\\
&=I_{i-1}\oplus
\begin{pmatrix}
a&a\\[1mm]
\dfrac{1-a}{a}&0
\end{pmatrix}
\oplus I_{n-i-1}.
\label{eq:burau-diagonal-conjugacy}
\end{align}
In the first line the upper-right entry is simply $r(1-a)=a$; it is written in this expanded form only to display the effect of diagonal conjugation.  Therefore, if $\beta_{1-a}$ denotes unreduced Burau with parameter $t=1-a$, the homogeneous matrix component of the present construction satisfies
\begin{equation}\label{eq:burau-equivalence}
\rho_{\mathrm{hom}}(\sigma_i)
=D_n\,\beta_{1-a}(\sigma_i)\,D_n^{-1}
\qquad(1\le i\le n-1).
\end{equation}
Thus the present family contains only the $b=a$ subfamily of the two-parameter 2022 construction literally, but its homogeneous specialization is representation-theoretically equivalent to generic unreduced Burau.  This comparison will be useful below: the Bigelow example collapses in the homogeneous Burau-equivalent regime, while a genuinely non-homogeneous parameter evaluation separates it.
\end{remark}

\section{Formal definitions and the semidirect-product representation}\label{sec:formal}

The local matrices are combined using a semidirect product so that the coefficients follow the moving strand labels.

\subsection{Coefficient field and local matrix algebra}\label{subsec:coefficient-field}

Let $k$ be a field and let $K=k(a_1,\ldots,a_n)$ be the rational function field in algebraically independent variables $a_1,\ldots,a_n$ (see, e.g., \cite{Lang2002}). By \eqref{eq:final-det}, each local matrix $\widehat{\sigma}_{i;xy}$ of \eqref{eq:final-hat} belongs to $GL_n(K)$.

\begin{definition}[Strand-labelled local matrix algebra]\label{def:An}
Let $\mathcal A_n\subseteq M_n(K)$ be the unital $K$-subalgebra generated by all $\widehat{\sigma}_{i;xy}$, where $1\le i\le n-1$ and $x,y\in[n]$ are distinct. Write $\mathcal A_n^\times$ for its unit group.
\end{definition}

Each local generator is a unit of $\mathcal A_n$: by Cayley--Hamilton, its inverse is a polynomial in that generator with coefficients in $K$.

\subsection{The label action of $S_n$}\label{subsec:label-action}

Every $P\in S_n$ acts on $K$ by $P\cdot a_x=a_{P(x)}$, extended to rational functions and entrywise to matrices. The formula for the local block gives
\begin{equation}\label{eq:local-equivariance}
P\cdot\widehat{\sigma}_{i;xy}=\widehat{\sigma}_{i;P(x)P(y)}.
\end{equation}
Consequently, $S_n$ acts on $\mathcal A_n^\times$ by group automorphisms.

\subsection{The semidirect-product (denoted by bullet product for simplicity)}\label{subsec:bullet}

\begin{definition}[The semidirect-product (bullet product)]\label{def:Gn}
On $\mathcal G_n=\mathcal A_n^\times\times S_n$, define the multiplication, denoted by $\bullet$, by
\begin{equation}\label{eq:bullet}
(X,P)\bullet(Y,Q)=\bigl(X(P\cdot Y),\,PQ\bigr).
\end{equation}
This defines the semidirect-product group
\begin{equation}\label{eq:Gn}
\mathcal G_n=\mathcal A_n^\times\rtimes S_n.
\end{equation}
The identity is $(I_n,\mathrm{id})$ and the inverse is
\begin{equation}\label{eq:bullet-inverse}
(X,P)^{-1}=\bigl(P^{-1}\cdot X^{-1},\,P^{-1}\bigr).
\end{equation}
\end{definition}

Thus $P$ relabels the coefficients of $Y$ before the matrix multiplication; it does not permute the matrix positions.

\subsection{The braid-group representation}\label{subsec:H-definition}

\begin{definition}[Generator assignment]\label{def:H-generators}
For $P_i=(i\;i+1)\in S_n$, prescribe
\begin{equation}\label{eq:H-generator}
H(\sigma_i)=\bigl(\widehat{\sigma}_{i;i,i+1},\,P_i\bigr),
\qquad 1\le i\le n-1.
\end{equation}
\end{definition}

For distant generators, the local blocks act on disjoint positions and their transpositions commute. For adjacent generators, the matrix equality is the general ordered-triple cubic identity of Proposition~\ref{prop:ordered-triple}; the corresponding permutation equality is the Coxeter relation. These facts give the following result.

\begin{theorem}[Strand-labelled local representation]\label{thm:representation}
For every $n\ge2$, the above assignment extends uniquely to a group homomorphism
\begin{equation}\label{eq:H-map}
\boxed{H:B_n\longrightarrow\mathcal G_n=\mathcal A_n^\times\rtimes S_n.}
\end{equation}
\end{theorem}
\begin{proof}
Straightforward.
\end{proof}

For negative generators, \eqref{eq:bullet-inverse} and \eqref{eq:local-equivariance} give
\begin{equation}\label{eq:H-inverse}
H(\sigma_i^{-1})
=\bigl(\widehat{\sigma}_{i;i+1,i}^{-1},\,P_i\bigr).
\end{equation}

\begin{remark}[Nature of the target]\label{rem:nature-target}
The representation takes values in the matrix--permutation group $\mathcal G_n$. Its matrix component uses ordinary $n\times n$ matrix multiplication, while the permutation component transports strand labels through successive crossings. The local coefficients therefore depend on the participating strands rather than solely on the Artin position.
\end{remark}

\section{Separation of Bigelow's Burau-kernel commutator}\label{sec:examples}

\paragraph{The local matrices in $B_5$.}
For the exact five-strand calculation below, the general rule \eqref{eq:final-hat} takes the following explicit forms at the first two Artin positions:
\begin{equation}\label{eq:sigma1xy-B5}
\widehat{\sigma}_{1;xy}
=
\begin{pmatrix}
a_x&a_y&0&0&0\\
\dfrac{1-a_x}{a_x}&0&0&0&0\\
0&0&1&0&0\\
0&0&0&1&0\\
0&0&0&0&1
\end{pmatrix},
\end{equation}
and
\begin{equation}\label{eq:sigma2xy-B5}
\widehat{\sigma}_{2;xy}
=
\begin{pmatrix}
1&0&0&0&0\\
0&a_x&a_y&0&0\\
0&\dfrac{1-a_x}{a_x}&0&0&0\\
0&0&0&1&0\\
0&0&0&0&1
\end{pmatrix}.
\end{equation}
The matrices at positions $3$ and $4$ have the same local block in the corresponding two rows and columns. Here $x,y$ denote the \emph{current strand labels}, not the row indices; they must be updated as crossings are read. Inverse letters are evaluated with the swapped-label inverse, as in \eqref{eq:H-inverse}. These conventions are used for Bigelow's braid word below.

We use Bigelow's explicit braid in $B_5$ \cite{Bigelow1999}.

\begin{definition}[Bigelow's commutator]\label{def:bigelow-commutator}
Define
\begin{equation}\label{eq:bigelow-psi1}
\psi_1=\sigma_3^{-1}\sigma_2\sigma_1^2\sigma_2\sigma_4^3\sigma_3\sigma_2,
\end{equation}
\begin{equation}\label{eq:bigelow-psi2}
\psi_2=\sigma_4^{-1}\sigma_3\sigma_2\sigma_1^{-2}\sigma_2\sigma_1^2\sigma_2^2\sigma_1\sigma_4^5,
\end{equation}
and
\begin{equation}\label{eq:bigelow-Sigma}
\Sigma_1=\psi_1^{-1}\sigma_4\psi_1,
\qquad
\Sigma_2=\psi_2^{-1}
\sigma_4\sigma_3\sigma_2\sigma_1^2\sigma_2\sigma_3\sigma_4
\psi_2.
\end{equation}
Let
\begin{equation}\label{eq:bigelow-C-main}
C=[\Sigma_1,\Sigma_2]
=\Sigma_1^{-1}\Sigma_2^{-1}\Sigma_1\Sigma_2.
\end{equation}
\end{definition}

\begin{remark}[Bigelow's kernel result]\label{rem:bigelow-kernel-result}
Bigelow proved that the commutator $C$ is nontrivial and belongs to the kernel of the Burau representation \cite{Bigelow1999}. It therefore provides a concrete test of the separation achieved by retaining individual strand parameters.
\end{remark}

\begin{theorem}[Non-homogeneous separation of Bigelow's commutator]\label{thm:bigelow-separation}
Over $k=\mathbb Q$, consider the exact specialization
\begin{equation}\label{eq:bigelow-values}
\boxed{
a_1=a_2=a_3=a_4=2,\qquad a_5=-1.
}
\end{equation}
If $H(C)=(X_C,P_C)$, then $P_C=\mathrm{id}$ and the evaluation of $X_C$ at these parameters is not $I_5$. In particular,
\[
\boxed{H(C)\neq (I_5,\mathrm{id})}
\]
over $\mathbb Q(a_1,\ldots,a_5)$. By contrast, the homogeneous specialization $a_1=\cdots=a_5=2$ is Burau-equivalent and sends $C$ to the identity.
\end{theorem}

\begin{proof}
Write $H(\beta)=(X_\beta,P_\beta)$ for the symbolic representation of Section~\ref{sec:formal}. The local blocks and their inverses involve only the parameters $a_r$ and denominators $a_r$ and $1-a_r$. Since products and strand relabellings preserve this form, every entry of $X_\beta$ belongs to the localization
\[
R=\mathbb Q[a_1,\ldots,a_5,\,a_1^{-1},\ldots,a_5^{-1},\,(1-a_1)^{-1},\ldots,(1-a_5)^{-1}].
\]
The values in \eqref{eq:bigelow-values} are neither $0$ nor $1$, so evaluation defines a ring homomorphism $\operatorname{ev}_{\alpha}:R\to\mathbb Q$. We abbreviate $\operatorname{ev}_{\alpha}(X_\beta)$ by $H(\beta)_{\mathrm{mat}}$ in the calculation below. This notation denotes an exact evaluation of the \emph{symbolic matrix component}, not a separate matrix representation of $B_5$: all strand-label permutations are applied before evaluation, according to the bullet product. Consequently, if two symbolic matrix components were equal, their evaluations would agree.

The braid words $\psi_1,\psi_2,\Sigma_1,\Sigma_2$ and $C$ are given in \eqref{eq:bigelow-psi1}--\eqref{eq:bigelow-C-main}. Figure~\ref{fig:bigelow-psis} depicts the two building braids.

\begin{figure}[H]
\centering
\begin{minipage}[t]{0.46\textwidth}
\centering
\resizebox{0.50\linewidth}{!}{%
\includegraphics{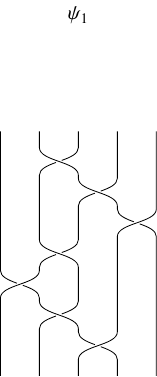}}
\end{minipage}
\hfill
\begin{minipage}[t]{0.46\textwidth}
\centering
\resizebox{0.50\linewidth}{!}{%
\includegraphics{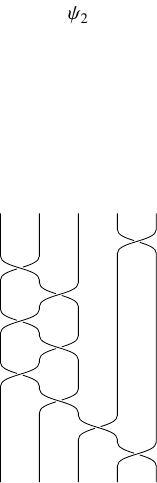}}
\end{minipage}
\caption{The two building braids in Bigelow's example, $\psi_1$ and $\psi_2$.  The algebraic words are exactly \eqref{eq:bigelow-psi1} and \eqref{eq:bigelow-psi2}.}
\label{fig:bigelow-psis}
\end{figure}

Both $\Sigma_1\Sigma_2$ and $\Sigma_2\Sigma_1$ have bottom strand order
\[
(s_1,s_4,s_3,s_2,s_5).
\]
Thus their permutation components agree. To determine whether the symbolic matrix components agree, it suffices to find one admissible evaluation at which they differ. The matrices below are obtained by applying the strand-label permutations during the symbolic bullet products and then evaluating the resulting rational functions at \eqref{eq:bigelow-values}.

\paragraph{Non-homogeneous evaluation.}
The choice \eqref{eq:bigelow-values} is non-homogeneous: strand $s_5$ has parameter $-1$ while the other four strands have parameter $2$.  Direct strand tracking and ordinary matrix multiplication give
\[
H(\psi_1)_{\mathrm{mat}}=
\begin{pmatrix}
3&8&-8&-8&0\\
-2&-6&5&4&0\\
0&1&0&0&0\\
\frac14&-\frac12&-\frac12&-\frac12&0\\
0&\frac38&0&0&-\frac12
\end{pmatrix},
\]
with bottom order $(s_1,s_5,s_2,s_4,s_3)$, and
\[
H(\psi_2)_{\mathrm{mat}}=
\begin{pmatrix}
6&7&12&0&0\\
-1&-1&0&0&0\\
-\frac32&-\frac32&-4&-6&-10\\
0&0&0&\frac54&2\\
\frac3{16}&\frac3{16}&\frac12&\frac58&1
\end{pmatrix},
\]
with bottom order $(s_2,s_1,s_5,s_4,s_3)$.

Continuing through the conjugations in \eqref{eq:bigelow-Sigma} gives
\begin{equation}\label{eq:bigelow-HS1}
H(\Sigma_1)_{\mathrm{mat}}=
\begin{pmatrix}
5&16&16&16&-64\\
-\frac34&-2&-3&-3&12\\
-\frac14&-1&0&-1&4\\
-\frac18&-\frac12&-\frac12&\frac12&2\\
-\frac1{32}&-\frac18&-\frac18&-\frac18&\frac32
\end{pmatrix},
\end{equation}
whose bottom order is $(s_1,s_4,s_3,s_2,s_5)$, while
\begin{equation}\label{eq:bigelow-HS2}
H(\Sigma_2)_{\mathrm{mat}}=
\begin{pmatrix}
3&-40&-96&-96&320\\
4&-5&-12&0&128\\
-\frac32&0&0&-4&-40\\
-\frac38&\frac92&\frac{43}{4}&10&-34\\
-\frac1{16}&1&\frac{19}{8}&\frac52&-8
\end{pmatrix},
\end{equation}
whose permutation component is the identity.  At the particular evaluation \eqref{eq:bigelow-values}, the transposition of labels $2$ and $4$ fixes all numerical parameter values, so the bullet products reduce numerically to the displayed ordinary products of \eqref{eq:bigelow-HS1} and \eqref{eq:bigelow-HS2}.  We obtain
\begin{equation}\label{eq:bigelow-prod12}
H(\Sigma_1\Sigma_2)_{\mathrm{mat}}=
\begin{pmatrix}
53&-272&-652&-544&2976\\
-\frac{43}{8}&\frac{77}{2}&\frac{369}{4}&84&-370\\
-\frac{37}{8}&\frac{29}{2}&\frac{139}{4}&24&-206\\
-\frac{31}{16}&\frac{47}{4}&\frac{225}{8}&24&-117\\
-\frac{29}{64}&\frac{45}{16}&\frac{215}{32}&6&-\frac{115}{4}
\end{pmatrix},
\end{equation}
and
\begin{equation}\label{eq:bigelow-prod21}
H(\Sigma_2\Sigma_1)_{\mathrm{mat}}=
\begin{pmatrix}
71&232&176&176&-768\\
\frac{91}{4}&70&63&75&-172\\
-\frac{23}{4}&-17&-17&-21&28\\
-\frac{65}{8}&-\frac{53}{2}&-\frac{81}{4}&-21&90\\
-\frac{55}{32}&-\frac{45}{8}&-\frac{17}{4}&-\frac{33}{8}&\frac{37}{2}
\end{pmatrix}.
\end{equation}
Their difference is the explicitly nonzero matrix
\begin{equation}\label{eq:bigelow-difference}
\boxed{
\begin{pmatrix}
-18&-504&-828&-720&3744\\
-\frac{225}{8}&-\frac{63}{2}&\frac{117}{4}&9&-198\\
\frac98&\frac{63}{2}&\frac{207}{4}&45&-234\\
\frac{99}{16}&\frac{153}{4}&\frac{387}{8}&45&-207\\
\frac{81}{64}&\frac{135}{16}&\frac{351}{32}&\frac{81}{8}&-\frac{189}{4}
\end{pmatrix}.
}
\end{equation}
In particular, the $(1,1)$ entries are $53$ and $71$, so the difference already has $(1,1)$ entry $-18$.

The exact difference in \eqref{eq:bigelow-difference} is nonzero. Because the evaluation homomorphism $\operatorname{ev}_{\alpha}$ preserves equality, the symbolic matrices in $M_5(\mathbb Q(a_1,\ldots,a_5))$ must also differ. Thus
\begin{equation}\label{eq:bigelow-separated}
\boxed{
H(\Sigma_1\Sigma_2)\neq H(\Sigma_2\Sigma_1).
}
\end{equation}
For the commutator $C$ defined in \eqref{eq:bigelow-C-main}, the permutation component of $H(C)$ is $\mathrm{id}$ and, at \eqref{eq:bigelow-values},
\begin{equation}\label{eq:bigelow-commutator}
H(C)_{\mathrm{mat}}=
\begin{pmatrix}
-\frac{5281}{2}&12222&29313&23616&-140040\\
\frac{22149}{32}&-\frac{25633}{8}&-\frac{122967}{16}&-6192&\frac{73431}{2}\\
\frac{5283}{32}&-\frac{6111}{8}&-\frac{29297}{16}&-1476&\frac{17505}{2}\\
\frac{3033}{64}&-\frac{3501}{16}&-\frac{16803}{32}&-422&\frac{10035}{4}\\
\frac{3483}{256}&-\frac{4023}{64}&-\frac{19305}{128}&-\frac{243}{2}&\frac{11545}{16}
\end{pmatrix}
\neq I_5.
\end{equation}
In particular, the $(1,1)$ entry of $H(C)_{\mathrm{mat}}$ is $-5281/2\neq1$, and hence the symbolic matrix component $X_C$ is not $I_5$. The commutator has trivial strand permutation, so $H(C)\neq(I_5,\mathrm{id})$, as asserted.

\paragraph{Homogeneous comparison.}
If all five parameters are set equal, the local family becomes the Burau-equivalent homogeneous family of Remark~\ref{rem:burau-connection}.  Direct calculation, for example at $a_1=\cdots=a_5=2$, gives
\[
H(\Sigma_1\Sigma_2)_{\mathrm{mat}}
=
H(\Sigma_2\Sigma_1)_{\mathrm{mat}},
\qquad
H(C)_{\mathrm{mat}}=I_5.
\]
This homogeneous comparison agrees with Bigelow's Burau-kernel construction \cite{Bigelow1999}, whereas the unequal-parameter representation distinguishes the commutator from the identity.
\end{proof}

\begin{remark}[Evidence toward faithfulness]\label{rem:bigelow-strength}
The separation in Theorem~\ref{thm:bigelow-separation} is a significant test of the strand-tracking representation: an explicit nontrivial braid annihilated by Burau is detected after the parameters are allowed to follow individual strands. The result provides strong evidence that the strand-tracking construction is a promising candidate for a faithful Burau-type representation. Establishing faithfulness requires a separate investigation beyond the scope of the present paper.
\end{remark}

\section{Conclusion}\label{sec:conclusion}

We have constructed a strand-labelled local representation of $B_n$ in which the crossing position specifies the location of a $2\times2$ block and the ordered strand labels determine its entries. Within the ansatz
\[
M_{xy}=\begin{pmatrix}a_x&a_y\\ b_{xy}&c_{xy}\end{pmatrix},
\]
the labelled cubic relations for arbitrary ordered triples determine $c_{xy}=0$ and $b_{xy}=(1-a_x)/a_x$ for $n\geq3$ in the nondegenerate branch. The same local family is used for $B_2$ by specialization. Keeping track of the moving labels leads to the semidirect-product representation $H:B_n\to\mathcal A_n^\times\rtimes S_n$. Its matrix component can be evaluated exactly at admissible parameter values after the strand-label relabellings have been performed.

In the special homogeneous specialization $a_x=a$ for every strand label $x$, the local core belongs to the $b=a$ subfamily of Case~1 in \cite{Pourkia2022}, and the resulting matrix representation is diagonally equivalent to unreduced Burau. The non-homogeneous specialization $a_1=a_2=a_3=a_4=2$, $a_5=-1$ produces a sharper distinction: it detects Bigelow's explicit nontrivial $B_5$ Burau-kernel commutator, with the nonidentity established by exact matrix calculations. This separation illustrates the role of moving strand labels in the representation and provides a concrete basis for further study of its kernel.

\end{document}